\documentclass[12pt,a4paper]{article}
\usepackage[cp1251]{inputenc}
\usepackage[english,russian]{babel}
\usepackage{amsmath}
\usepackage{amssymb}
\usepackage{amsfonts}
\usepackage{hyperref}

\newcounter{paragr}

\begin{document}

\binoppenalty=10000 \relpenalty=10000

\renewcommand{\refname}{References}
\renewcommand{\contentsname}{Contents}

\begin{center}
{\huge Viscous compressible homogeneous multi-fluids with multiple velocities: barotropic existence theory}
\end{center}

\medskip

\begin{center}
{\large Alexander Mamontov\footnote{The author is partially supported by the Russian Foundation for Basic Research (project 15--01--08275).},\quad Dmitriy Prokudin}
\end{center}

\medskip

\begin{center}
{\large October 18, 2016}
\end{center}

\medskip

\begin{center}
{
Lavrentyev Institute of Hydrodynamics, \\ Siberian Branch of the Russian Academy of Sciences\\
pr. Lavrent'eva 15, Novosibirsk 630090, Russia}
\end{center}

\medskip

\begin{center}
{\bfseries Abstract}
\end{center}


\begin{center}
\begin{minipage}{110mm}
 We consider the model of viscous compressible homogeneous multi-fluids with multiple velocities. We review different formulations of the model and the existence results for boundary value problems. We analyze crucial mathematical difficulties which arise during the proof of the existence theorems.
\end{minipage}
\end{center}

\bigskip


\newpage

\tableofcontents

\bigskip

\section{Introduction}
\stepcounter{paragr}\setcounter{equation}{0}\medskip

\noindent\indent Since the publication of our preceding review \cite{mamprok.france} 3 years ago an essential de\-ve\-lop\-ment has been achieved in the solvability theory for the equations of motion of multi-fluids, at least for one of the directions described in \cite{mamprok.france}. In the article, we give a description of the corresponding class of problems and give the formulations of the obtained results.

The work is written after the report made at the 8th International Scientific School-Conference of Young Scientists ``Theory and Numerical Methods of Solving Inverse and Ill-posed Problems'' (Novosibirsk, 1-7 September, 2016), the presentation of the report is available at the link \cite{mamprok.conf}.

\section{The model of multi-fluids with multiple velocities}
\stepcounter{paragr}\setcounter{equation}{0}\medskip

We do not go into details here explaining the formulation of the equations of viscous compressible homogeneous multi-fluids with multiple velocities, because it is de\-s\-cribed in \cite{mamprok.france} as well as in the papers which contain the corresponding solvability results (the references are given below). Let us immediately write down the final formulation of the corresponding system of PDEs:
\begin{equation}\label{mamprok.continit}
\frac{\partial\rho_{i}}{\partial t}+{\rm div\,}(\rho_{i}\boldsymbol{u}_{i})=0,\quad i=1, \ldots, N,
\end{equation}
\begin{equation}\label{mamprok.mominit}
\frac{\partial(\rho_{i}\boldsymbol{u}_{i})}{\partial t}+{\rm div\,}(\rho_{i}\boldsymbol{u}_{i}\otimes\boldsymbol{u}_{i})+\nabla p_{i}
={\rm div\,}{\mathbb S}_{i}+\boldsymbol{J}_{i}+\rho_{i}\boldsymbol{f}_{i},\quad i=1, \ldots, N.
\end{equation}
Here $\rho_{i}$ is the density of the $i$-th component (constituent) of the multi-fluid (totally $N$ components are present), $\boldsymbol{u}_{i}$ is the velocity field, $p_{i}$ is the pressure, ${\mathbb S}_{i}$ is the viscous stress tensor, the vectors $\displaystyle\boldsymbol{J}_{i}=\sum\limits_{j=1}^{N} a_{ij}(\boldsymbol{u}_{j}-\boldsymbol{u}_{i})$ are responsible for the intensity of the momentum exchange between the constituents of the multi-fluid, and the vectors~$\boldsymbol{f}_{i}$ are known fields of external body forces. The viscous stress tensors ${\mathbb S}_{i}$ are defined by the equalities
\begin{equation}\label{mamprok.tensorvisc}
{\mathbb S}_{i}=\sum\limits_{j=1}^{N}\widehat{{\mathbb S}}_{ij},\quad \widehat{{\mathbb S}}_{ij}=\Big(2\mu_{ij}{\mathbb D}(\boldsymbol{u}_{j})+\lambda_{ij}({\rm div\,}
\boldsymbol{u}_{j}){\mathbb I}\Big),\qquad  i, j=1, \ldots, N,
\end{equation}
where ${\mathbb D}(\boldsymbol{v})=((\nabla\otimes\boldsymbol{v})+(\nabla\otimes\boldsymbol{v})^{*})/2$ is the rate of deformation tensor of the vector field~$\boldsymbol{v}$,
${\mathbb I}$ is the identity tensor, and the viscosity coefficients compose the matrices
\begin{equation}\label{mamprok.uslovie_visc}
\textbf{M}=\{\mu_{ij}\}_{i, j = 1}^{N}>0,\quad \textbf{H}=\boldsymbol{\boldsymbol{\Lambda}}+\frac{2}{3}\textbf{M}\geqslant 0,\quad \boldsymbol{\Lambda}=\{\lambda_{ij}\}_{i,j= 1}^{N}.
\end{equation}

Generally speaking, the pressures $p_{i}$ depend not only on $\rho_{i}$, but also on other thermodynamical variables (say, temperatures), so that the system \eqref{mamprok.continit},
\eqref{mamprok.mominit} must be complemented by the equations for these variables. The resulting system describes motions of heat-conductive multi-fluids and, strictly speaking, it is this system that is consistent thermodynamically. This system was the subject of some existence results \cite{mamprok.smj12}, \cite{mamprok.izvran} (oddly enough, the theory for the barotropic system \eqref{mamprok.continit}, \eqref{mamprok.mominit} is less developed). However, the formulated model presents essential difficulties which lead to restrictions on the structure of the total viscosity matrix $\textbf{N}=\boldsymbol{\Lambda}+2\textbf{M}=\{\nu_{ij}\}_{i, j = 1}^{N}$. From the mathematical point of view, the only interesting case is that of non-diagonal viscosity matrices (which corresponds to the presence of the viscous friction between the constituents of the multi-fluid), when the inter\-penetration of the higher-order terms takes place in the equations \eqref{mamprok.mominit}, and the theory does not follow from the corresponding theory of mono-fluids (i.~e. from the usual NS or NSF theory). In the framework of the complete multiple-velocities statement, posed in \eqref{mamprok.continit}, \eqref{mamprok.mominit}, global existence results are obtained only for triangular matrices $\textbf{N}$ (the structure of the matrices $\boldsymbol{\Lambda}$ and $\textbf{M}$ is irrelevant). On the other hand, physical arguments lead to the symmetry of the viscosity matrices. Hence, the results completely reasonable in the physical sense are obtained only in the less interesting case of diagonal matrices $\textbf{N}$. Thus, further development of the mathematical theory of multi-fluids with multiple velocities (which would agree with the physical sense of solutions) is limited by essential mathematical difficulties, and maybe it requires correcting the model itself. More detailed description of the problem may be found in \cite{mamprok.france}.

Recently we discovered a variation of the model \eqref{mamprok.continit}, \eqref{mamprok.mominit}, in which the difficulties described above can be overcome, and all restrictions on the structure of the viscosity matrices can be removed, except physically necessary conditions \eqref{mamprok.uslovie_visc}. The corresponding results refer to the barotropic case, i.~e. to the (modified) system \eqref{mamprok.continit}, \eqref{mamprok.mominit}, in which every $p_{i}$ depends only on $\rho_{i}$. The exact formulation of the model is presented in the next Section.

\section{Model with common material derivative operator and common pressure}
\stepcounter{paragr}\setcounter{equation}{0}\medskip

Let us make the following corrections and accept the following additional as\-sump\-tions for the system \eqref{mamprok.continit}, \eqref{mamprok.mominit}:
\begin{itemize}
  \item The velocities $\boldsymbol{u}_{i}$ of the constituents in the convective terms, more exactly, in the material derivative operators
 $\displaystyle \frac{d}{dt}=\frac{\partial}{\partial t}+\boldsymbol{u}_{i}\cdot\nabla$, are replaced by the average velocity $\boldsymbol{v}$ of the multi-fluid.
  \item The pressures in all constituents are equal to each other: $p_i=p$, $i=1, \ldots, N$, and this common pressure depends only on the total density $\displaystyle\rho=\sum\limits_{i=1}^N \rho_{i}$ of the multi-fluid.
\end{itemize}
This leads to the following system of equations
\begin{equation}\label{mamprok.contfin}
\frac{\partial\rho_{i}}{\partial t}+{\rm div\,}(\rho_{i}\boldsymbol{v})=0,\quad i=1, \ldots, N,
\end{equation}
\begin{equation}\label{mamprok.momfin}
\frac{\partial(\rho_{i}\boldsymbol{u}_{i})}{\partial t}+{\rm div\,}(\rho_{i}\boldsymbol{v}\otimes\boldsymbol{u}_{i})+\nabla p(\rho)
={\rm div\,}{\mathbb S}_{i}+\rho_{i}\boldsymbol{f}_{i},\quad i=1, \ldots, N,
\end{equation}
instead of \eqref{mamprok.continit}, \eqref{mamprok.mominit}. We have also eliminated the terms $\boldsymbol{J}_{i}$, which is in full accord with the idea of closeness of $\boldsymbol{u}_{i}$ to each other. That is what the first of the above assumptions is based on (though, retention of $\boldsymbol{J}_{i}$ would not create any mathematical difficulties). The second assumption is quite popular in the theory of mixtures (see more details e.~g. in \cite{mamprok.semi1}, \cite{mamprok.izvrannew}). The relations \eqref{mamprok.tensorvisc} and \eqref{mamprok.uslovie_visc} hold true, no additional restrictions being necessary for the viscosity matrices apart from \eqref{mamprok.uslovie_visc}. In order to make the system \eqref{mamprok.contfin}, \eqref{mamprok.momfin} closed, it is left to specify the representation of the average velocity $\boldsymbol{v}$ as a function of the velocities $\boldsymbol{u}_{i}$ of the constituents. For the sake of simplicity we assume the relation
$\displaystyle\boldsymbol{v}=\frac{1}{N}\sum\limits_{i=1}^N \boldsymbol{u}_{i}$, but it is possible to consider more general and physically adequate versions in which the relative concentrations of the constituents would be taken into account.

It is possible to rewrite the system \eqref{mamprok.contfin}, \eqref{mamprok.momfin} in terms of concentrations $\alpha_i=\rho_i/\rho$, i.~e. to replace the equations
\eqref{mamprok.contfin} by their sum (the equation for~$\rho$) and transport equations for all concentrations $\alpha_i$ except one (say, $\alpha_N$). It is interesting to note that these two versions of writing down are completely equivalent in the unsteady case, whereas in the steady case there arise additional difficulties concerning the normalization of the  constituents masses and the uniqueness of solution. It is possible to observe this problem in more detail in \cite{mamprok.semi2}, where the steady system \eqref{mamprok.contfin}, \eqref{mamprok.momfin} is studied exactly in terms of concentrations. We are not aware of a well-posed formulation of the steady problem \eqref{mamprok.contfin}, \eqref{mamprok.momfin}, i.~e. such formulation, in which the unique definition of the densities $\rho_i$ in the steady state $\boldsymbol{u}_{i}=0$ would be provided. The problem is in the difference between the unsteady case, where the mass of each constituent and its distribution in the flow domain are uniquely defined by the initial data, and the steady state (which is understood as the limit state after the stabilization of an unsteady flow), in which it is not clear how the unsteady distributions of the densities are inherited. Formally, the stabilization leads to the relations $\boldsymbol{v}\cdot\nabla\alpha_i=0$ for the concentrations $\alpha_i$, which would be the source for their determination, but these relations are inconvenient as regards mathematics.

We can formulate the initial boundary value problem for the system \eqref{mamprok.contfin}, \eqref{mamprok.momfin} to describe unsteady flows in a bounded domain $\Omega\subset{\mathbb R}^3$ and prove the global existence theorems of weak solutions in the cylinder $Q_T=\Omega\times(0,T)$ for arbitrary $T>0$. In the polytropic case $p=K\rho^\gamma$, $\gamma>3/2$, it is made in \cite{mamprok.semi1}, and this result is generalized in \cite{mamprok.izvrannew} for the case of quite arbitrary dependence $p(\rho)$.

We can formulate the boundary value problem for the steady system
\begin{equation}\label{mamprok.contfinstat}
{\rm div\,}(\rho_{i}\boldsymbol{v})=0,\quad i=1, \ldots, N,
\end{equation}
\begin{equation}\label{mamprok.momfinstat}
{\rm div\,}(\rho_{i}\boldsymbol{v}\otimes\boldsymbol{u}_{i})+\nabla p(\rho)
={\rm div\,}{\mathbb S}_{i}+\rho_{i}\boldsymbol{f}_{i},\quad i=1, \ldots, N
\end{equation}
to describe steady flows in a bounded domain $\Omega\subset{\mathbb R}^3$ and prove the existence theorems of weak solutions. In the polytropic case $p=K\rho^\gamma$, $\gamma>3/2$, it is made in \cite{mamprok.semi2}, and this result is generalized in \cite{mamprok.smz1}, \cite{mamprok.smz2} for the case of quite arbitrary dependence~$p(\rho)$. Though, as it is mentioned above, in \cite{mamprok.semi2} rather a different statement of the problem is considered.

Let us note that the assumptions written down at the beginning of the Section were first used in \cite{mamprok.prokkraj} in the model case of equal densities of the constituents.

Below we give short comments on the main difficulties which we come across during the proof of solvability of boundary value problems for the systems \eqref{mamprok.continit}, \eqref{mamprok.mominit}, or \eqref{mamprok.contfin}, \eqref{mamprok.momfin}, or \eqref{mamprok.contfinstat}, \eqref{mamprok.momfinstat}, and on the difference between these difficulties and the corresponding aspects of the mono-fluid theory.

\section{Approximate solutions and compactness}
\stepcounter{paragr}\setcounter{equation}{0}\medskip

The proof of solvability for the boundary value problems mentioned in the preceding Section is made via the scheme similar to that in the mono-fluid theory, described in \cite{mamprok.feir09}, \cite{mamprok.novstrs04}. We do not reproduce here this scheme in full (it is partially described in the next Section), but just say that the solutions to the corresponding approximate problems are constructed, for which we obtain estimates uniform with respect to the approximating parameters, and, basing on these es\-ti\-mates, we deduce the weak convergence of the approximate solutions to some limits, which are going to be the solutions to the original problem. Thus, to within the terms inserted into the approximate equations, which as usual do not constitute difficulties greater than the terms permanently present in the (original) equations, the problem is reduced to the limit in the equations which look no different from the original ones.

In other words, the problem is reduced to the proof of the following conditional result. We assume that there is a sequence of solutions to the original problem (and they satisfy a priori estimates which provide weak convergence to some limits, after the necessary selection of a subsequence). It is required to prove that these limits themselves are a solution to the same problem. This is exactly what is called the compactness property of the set of solutions (relative compactness follows from the estimates, and the problem essence is in the closedness). If this kind of result is proved then the proof of the existence theorem is reduced to a skilled construction of approximate solutions in such a manner that would prevent appearance of new essential difficulties except for those that were overcome during the proof of the compactness of the solution set.

Such technique of constructing approximate solutions is well developed in the theory of viscous compressible mono-fluids, but it is not transportable completely and/or automatically to the multi-fluid case. The same concerns the crucial point in the proof of the compactness of the solutions set, that is the effective viscous flux technique. Thus, our goal would be to describe two kinds of difficulties and their solutions:
\begin{itemize}
  \item peculiarities which arise in the analysis of the effective viscous fluxes of the multi-fluid constituents,
  \item new effects which arise in the approximate equations of multi-fluid motions.
\end{itemize}
The corresponding comments are given in the two following Sections.

\section{Effective viscous fluxes and mixed products in the multi-fluid equations}

Similarly to mono-fluids, let us call the values $\displaystyle F_i=p_i-\sum\limits_{k=1}^N \nu_{ik}{\rm div}\boldsymbol{u}_{k}$ effective viscous fluxes of the multi-fluid constituents. In contrast with mono-fluids, in which this flux is unique, here we have $N$ fluxes, separate ones for each constituent. By way of illustration let us first deal with these values in the case of the steady version of the system \eqref{mamprok.continit}, \eqref{mamprok.mominit} (i.~e. the version in which the terms containing the derivatives with respect to $t$ are rejected) and its modified version \eqref{mamprok.contfinstat}, \eqref{mamprok.momfinstat}. The unsteady version adds some peculiarities which would be commented on in the next Section. We omit technical details of justifying some operations, implying them valid. So then, let us consider an assumed solution to the system
\begin{equation}\label{mamprok.contcommon}
{\rm div\,}(\rho_{i}\boldsymbol{w}_{i})=0,\quad i=1, \ldots, N,
\end{equation}
\begin{equation}\label{mamprok.momcommon}
{\rm div\,}(\rho_{i}\boldsymbol{w}_{i}\otimes\boldsymbol{u}_{i})+\nabla p_{i}
={\rm div\,}{\mathbb S}_{i}+\boldsymbol{z}_{i}+\rho_{i}\boldsymbol{f}_{i},\quad i=1, \ldots, N,
\end{equation}
in which $\boldsymbol{w}_{i}=\boldsymbol{u}_{i}$, $\boldsymbol{z}_{i}=\boldsymbol{J}_{i}$ for the original model, and $\boldsymbol{w}_{i}=\boldsymbol{v}$, $\boldsymbol{z}_{i}=0$, $p_i=p$ for the modified model. The boundary conditions are accepted in the form $\boldsymbol{u}_{i}|_{\partial\Omega}=0$. It is necessary to justify the weak limit in the written system, and the unique difficulty is contained in the terms $p_i$: after the limit they turn into $\overline{p}_i$ (here and below the bar denotes the weak limit), and it is necessary to justify the equalities $\overline{p}_i={p}_i$.

For all smooth functions $\tau$ and for all $i,j=1,\ldots,N$ we have the identities
$${\rm div}\left({\mathbb S}_{i}(\Delta^{-1}\rho_{j})\nabla\tau+{\mathbb S}_{i}\tau\nabla\Delta^{-1}\rho_j-{\rm div}{\mathbb S}_{i}\cdot\tau\Delta^{-1}\rho_{j}\right)=$$
$$={\mathbb S}_{i}:(\nabla\otimes[(\Delta^{-1}\rho_{j})\nabla\tau])+{\mathbb S}_{i}:(\nabla\otimes[\tau\nabla\Delta^{-1}\rho_j])-
({\rm div}{\rm div}{\mathbb S}_{i})\tau\Delta^{-1}\rho_{j},$$
which turn into the relations
\begin{equation}\label{mamprok.ident}
{\mathbb S}_{i}:(\nabla\otimes[(\Delta^{-1}\rho_{j})\nabla\tau])+{\mathbb S}_{i}:(\nabla\otimes[\tau\nabla\Delta^{-1}\rho_j])\overset{\Omega}{\sim}
({\rm div}{\rm div}{\mathbb S}_{i})\tau\Delta^{-1}\rho_{j},
\end{equation}
provided that $\tau\in C^{\infty}_{0}(\Omega)$, where $\overset{\Omega}{\sim}$ stands for equalities to within terms dis\-ap\-pearing after the integration over $\Omega$. Let us note that
$\displaystyle {\rm div}{\rm div}{\mathbb S}_{i}=\sum\limits_{k=1}^N \nu_{ik}\Delta{\rm div}\boldsymbol{u}_{k}$. Multiplying~\eqref{mamprok.momcommon} by
$\tau\nabla\Delta^{-1}\rho_j$ and taking into account \eqref{mamprok.ident}, we obtain the relations
\begin{equation}\label{mamprok.eff}
\tau\rho_j F_i\overset{\Omega,\tau}{\sim}-\tau(\rho_{i}\boldsymbol{w}_{i}\otimes\boldsymbol{u}_{i}):(\nabla\otimes\nabla\Delta^{-1}\rho_j)-
\tau(\boldsymbol{z}_{i}+\rho_{i}\boldsymbol{f}_{i})\nabla\Delta^{-1}\rho_j,
\end{equation}
where $\overset{\Omega,\tau}{\sim}$ stands for equalities to within terms disappearing after the integration over $\Omega$ and/or containing the derivatives of $\tau$ and hence being lower-order with respect to the unknowns (and hence constituting no difficulties during the limit, i.~e. disappearing in communicative relations). Let us accept below $\tau=1$ (in fact, this means additional work with lower-order terms, which does not present difficulties, and further limit $\tau\to 1$).

Let us introduce into consideration the operator ${\rm Comm}$ which acts
as\footnote{The operator ${\rm Comm}$ maps scalar functions into symmetric second rank tensors. Acting of this operator on non-scalar arguments implies the sum via some indices, and there is no necessity to specify these indices because of the symmetry. Thus, if ${\rm Comm}$ acts on vector-valued functions, the operator
$\nabla\otimes\nabla\Delta^{-1}$ is understood as $\nabla\Delta^{-1}{\rm div}$.}
$${\rm Comm}(a,b)=(\nabla\otimes\nabla\Delta^{-1}a)b-a(\nabla\otimes\nabla\Delta^{-1}b).$$
Then we obtain by definition
\begin{equation}\label{mamprok.eff11}
{\rm Comm}(\rho_{i}\boldsymbol{u}_{i},\rho_j)=(\nabla\Delta^{-1}{\rm div\,}(\rho_{i}\boldsymbol{u}_{i}))\rho_j-(\nabla\otimes\nabla\Delta^{-1}\rho_j)\rho_{i}\boldsymbol{u}_{i},
\end{equation}
and hence, since the operator $\nabla\otimes\nabla\Delta^{-1}$ is self-conjugate, we deduce
\begin{equation}\label{mamprok.eff12}
\boldsymbol{w}_{i}\cdot{\rm Comm}(\rho_{i}\boldsymbol{u}_{i},\rho_j)\overset{\Omega}{\sim}\rho_{i}\boldsymbol{u}_{i}\cdot \nabla\Delta^{-1}{\rm div\,}(\rho_{j}\boldsymbol{w}_{i})-
(\rho_{i}\boldsymbol{w}_{i}\otimes\boldsymbol{u}_{i}):(\nabla\otimes\nabla\Delta^{-1}\rho_j).
\end{equation}

Now, due to \eqref{mamprok.contcommon} and \eqref{mamprok.eff11}, the relations \eqref{mamprok.eff} take the form
\begin{equation}\label{mamprok.eff1}
\rho_j F_i\overset{\Omega,\tau}{\sim}\boldsymbol{w}_{i}\cdot{\rm Comm}(\rho_{i}\boldsymbol{u}_{i},\rho_j)-(\boldsymbol{z}_{i}+\rho_{i}\boldsymbol{f}_{i})\nabla\Delta^{-1}\rho_j.
\end{equation}
The right-hand sides of these relations do not contain products of weakly convergent values, except the bilinear form ${\rm Comm}$, which is known to survive the weak limit. Hence \eqref{mamprok.eff1} lead to the communicative relations
\begin{equation}\label{mamprok.eff2}
\rho_j \overline{F}_i\overset{\Omega}{\sim}\overline{\rho_j F_i},\quad i,j=1, \ldots, N.
\end{equation}
In order to prove the relations $\overline{p}_i=p_i$, it suffices to exclude the velocities from \eqref{mamprok.eff2}, and after that we can use standard arguments based on the monotonicity of the pressure(s) as (a) function(s) of the density(-ies). For instance, it would be sufficient to prove the relations
\begin{equation}\label{mamprok.eff3}
\sum\limits_{k=1}^N \nu_{ik}\overline{\rho_j{\rm div}\boldsymbol{u}_{k}}\overset{\Omega}{\sim}\sum\limits_{k=1}^N \nu_{ik}\rho_j{\rm div}\boldsymbol{u}_{k}
\end{equation}
at least for certain $i$, $j$ (in the case of the original model \eqref{mamprok.continit}, \eqref{mamprok.mominit} it would be sufficient to consider $i=j$). However, generally, the proof of \eqref{mamprok.eff3} is inconvenient. The renormalization of the equations \eqref{mamprok.contcommon} leads to the relations
$\rho_{i}{\rm div\,}\boldsymbol{w}_{i}\overset{\Omega}{\sim}0$, $i=1, \ldots, N$, from which it follows that
\begin{equation}\label{mamprok.eff4}
\overline{\rho_{i}{\rm div\,}\boldsymbol{w}_{i}}\overset{\Omega}{\sim}\rho_{i}{\rm div\,}\boldsymbol{w}_{i}.
\end{equation}
The following variants of further argument are possible:
\begin{itemize}
  \item Let us consider the original model \eqref{mamprok.continit}, \eqref{mamprok.mominit}, then \eqref{mamprok.eff4} look as
  $\overline{\rho_{i}{\rm div\,}\boldsymbol{u}_{i}}\overset{\Omega}{\sim}\rho_{i}{\rm div\,}\boldsymbol{u}_{i}$. If the matrix $\textbf{N}$ is diagonal, then we obtain immediately
  \eqref{mamprok.eff3} for all $i=j$, and hence $\rho_i \overline{p}_i\overset{\Omega}{\sim}\overline{\rho_i p_i}$ for all $i$, which leads to the strong convergence of the densities. If the matrix $\textbf{N}$ is triangular (say, upper triangular), then we first obtain the strong convergence of only one density (say, $\rho_N$), but further we use similar arguments to consider all densities step by step. This scheme is used in \cite{mamprok.smj12} and \cite{mamprok.izvran}. For matrices~$\textbf{N}$ of a general form there appear difficulties which are noticeable in the argument above, they were also discussed in \cite{mamprok.france}.
\item Let us consider the modified model \eqref{mamprok.contfinstat}, \eqref{mamprok.momfinstat}. Summing \eqref{mamprok.eff4} over $i=1, \ldots, N$, we obtain
$\overline{\rho{\rm div\,}\boldsymbol{v}}\overset{\Omega}{\sim}\rho{\rm div\,}\boldsymbol{v}$, and after the summation of \eqref{mamprok.eff2} over $j=1, \ldots, N$ we come to the relations $\rho\overline{F}_i\overset{\Omega}{\sim}\overline{\rho F_i}$, from which, due to the preceding relation, the velocities can be excluded, and finally
$\rho \overline{p}\overset{\Omega}{\sim}\overline{\rho p}$, which leads to the strong convergence of the density. Such scheme is used in
\cite{mamprok.semi1}, \cite{mamprok.semi2}, \cite{mamprok.izvrannew}, \cite{mamprok.smz1}, \cite{mamprok.smz2}, \cite{mamprok.prokkraj}.
\end{itemize}

The arguments above show a fundamental difference between multi-fluids and mono-fluids, namely, the presence of the mixed products of the form $\rho_{i}{\rm div\,}\boldsymbol{u}_{j}$, which are impossible to analyze with the use of the continuity equations \eqref{mamprok.contcommon}. Hence, in order to apply the effective viscous flux technique, additional tricks or assumptions for the model are required.

One more peculiarity appears in the proof of the compactness for sufficiently small $\gamma$ (i.~e. the rate of growth of the pressure(s) close to $\rho^{3/2}$ or $\rho_i^{3/2}$ respectively). In that case not the functions $\rho$ and $\rho_i$ are used in \eqref{mamprok.eff4}, but their cut-offs $T_r(\rho)$ or $T_r(\rho_i)$. The cut-off function $T_r(s)=s\chi_{s<r}(s)+r\chi_{s>r}(s)$ is nonlinear, hence \eqref{mamprok.eff4} with $T_r(\rho)$ could not be derived via the summation over $i$ (however, such relation could be derived via direct renormalization of the equation for the total density). Correspondingly, instead of \eqref{mamprok.eff2}, one should derive the relations in which
the corresponding cut-off ($T_r(\rho_i)$ in the original model or $T_r(\rho)$ in the modified one) is taken instead of $\rho_j$, since the summation over $j$ do not lead to the desired result either.

\section{Unsteady problems}
\stepcounter{paragr}\setcounter{equation}{0}\medskip

Let us comment on the peculiarities which arise in the proof of the compactness of the solutions set for the system
\begin{equation}\label{mamprok.contcommonnonst}
\frac{\partial\rho_{i}}{\partial t}+{\rm div\,}(\rho_{i}\boldsymbol{w}_{i})=0,\quad i=1, \ldots, N,
\end{equation}
\begin{equation}\label{mamprok.momcommonnonst}
\frac{\partial(\rho_{i}\boldsymbol{u}_{i})}{\partial t}+{\rm div\,}(\rho_{i}\boldsymbol{w}_{i}\otimes\boldsymbol{u}_{i})+\nabla p_{i}
={\rm div\,}{\mathbb S}_{i}+\boldsymbol{z}_{i}+\rho_{i}\boldsymbol{f}_{i},\quad i=1, \ldots, N
\end{equation}
in comparison with \eqref{mamprok.contcommon}, \eqref{mamprok.momcommon}. The sense of the values $\boldsymbol{w}_{i}$, $\boldsymbol{z}_{i}$ and $p_i$, as well as the boundary conditions, remain the same.

Let us multiply \eqref{mamprok.momcommonnonst} by $\psi\tau\nabla\Delta^{-1}\rho_j$, where $\psi\in C^{\infty}_{0}(0,T)$, and use \eqref{mamprok.ident} again. This leads us to the relations similar to \eqref{mamprok.eff}, but with the following differences:
\begin{itemize}
  \item the factor $\tau$ is substituted by $\psi\tau$,
  \item $\overset{\Omega,\tau}{\sim}$ is substituted by $\overset{Q_T,\tau,\psi}{\sim}$, which means equalities to within terms dis\-ap\-pearing after the integration over $Q_T$ and/or containing the derivatives of $\tau$ and/or $\psi$ (with a similar effect),
  \item the following term is added to the left-hand side (and it may be trans\-formed using \eqref{mamprok.contcommonnonst} and the fact that the operator $\nabla\otimes\nabla\Delta^{-1}$ is self-conjugate):
$$-\frac{\partial(\rho_{i}\boldsymbol{u}_{i})}{\partial t}\psi\tau\cdot\nabla\Delta^{-1}\rho_j\overset{Q_T,\tau,\psi}{\sim}\rho_{i}\boldsymbol{u}_{i}\psi\tau\cdot
\nabla\Delta^{-1}\frac{\partial\rho_{j}}{\partial t}\overset{Q_T,\tau,\psi}{\sim}
$$
$$
-\psi\tau\rho_j\boldsymbol{w}_{j}\cdot
\nabla\Delta^{-1}{\rm div}(\rho_{i}\boldsymbol{u}_{i}).$$
\end{itemize}

On the other hand, instead of \eqref{mamprok.eff12}, we use the identities
\begin{equation}\label{mamprok.eff13}
\boldsymbol{w}_{i}\cdot{\rm Comm}(\rho_{i}\boldsymbol{u}_{i},\rho_j)\overset{\Omega}{\sim}
\rho_{j}\boldsymbol{w}_{i}\cdot\nabla\Delta^{-1}{\rm div}(\rho_{i}\boldsymbol{u}_{i})-
(\rho_{i}\boldsymbol{w}_{i}\otimes\boldsymbol{u}_{i}):(\nabla\otimes\nabla\Delta^{-1}\rho_j),
\end{equation}
in other words, we use \eqref{mamprok.eff11} in the original form, without the transfer of the operator $\nabla\otimes\nabla\Delta^{-1}$ from $\rho_{j}\boldsymbol{w}_{i}$ to $\rho_{i}\boldsymbol{u}_{i}$.

As a result, instead of \eqref{mamprok.eff1}, we obtain the relations (we set again $\tau=1$ and similarly $\psi=1$)
\begin{equation}\label{mamprok.eff5}\begin{array}{c}
\rho_j F_i\overset{Q_T,\tau,\psi}{\sim}\boldsymbol{w}_{i}\cdot{\rm Comm}(\rho_{i}\boldsymbol{u}_{i},\rho_j)-(\boldsymbol{z}_{i}+\rho_{i}\boldsymbol{f}_{i})\nabla\Delta^{-1}\rho_j+
\\ \\ +\rho_{j}(\boldsymbol{w}_{j}-\boldsymbol{w}_{i})\cdot\nabla\Delta^{-1}{\rm div}(\rho_{i}\boldsymbol{u}_{i}),
\end{array}\end{equation}
The last term disappears for the modified model, and \eqref{mamprok.eff5} turns into \eqref{mamprok.eff1}. For the original model, we take only $j=i$, that gives
\eqref{mamprok.eff2} for all $i=j$, and it suffices.

Further argument (during the proof of the compactness of the solutions set) differs from that described in Section 4 only in typical difficulties which distinguish unsteady problems from the steady ones, and which are not a subject of our paper. In other words, the text of Section 4 starting from the relations \eqref{mamprok.eff2} and up to the end remains the same for the unsteady case.

During the constructing of the approximate solutions, there arise additional difficulties related to the estimates and convergence of the mixed products of the form $\rho_{i}\boldsymbol{u}_{j}$, which are impossible to analyze via the momentum equations \eqref{mamprok.momcommon}. This leads to the necessity of estimation of the ratios $\rho_i/\rho_j$, see the details e.~g. in \cite{mamprok.semi1}, \cite{mamprok.izvrannew}.

\section{One-dimensional problems}
\stepcounter{paragr}\setcounter{equation}{0}\medskip

In spite of the essential progress in the multidimensional viscous gas theory, the one-dimensional theory with its apex in 1970-80s, did not lose its relevance up to now. This is due to at least two factors. Firstly, the multidimensional existence theorems concern only weak solutions, whose regularity is not sufficient even for the uniqueness; smoothness increase is hindered by serious obstacles. Secondly, the difficulty of multidimensional problems eclipses the study of many qualitative properties of solutions, as well as related problems including modeling; this gives researchers the right to consider the corresponding questions first in the one-di\-men\-sional case. Hence, whereas in the solvability theory for the main boundary value problems of the viscous gas there was a shift of emphasis from the one-dimensional case to the multi\-dimensional one already two decades ago, in many other domains of the theory, the one-dimensional problems stay at the forefront.

Everything mentioned above concerns multi-fluids as well. As it is clear from the paper, the solvability theory for multi-fluids is now comparable with that for mono-fluids. At the same time some problems appeared concerning multi-fluids specifically and distinguishing them from mono-fluids. A logical thing to do would be to study the equations of one-dimensional motions of homogeneous viscous compressible multi-fluids with multiple velocities starting with the existence and uniqueness theorems for boundary value problems with a view to using it as a basis for studying the statement of the problems and studying properties of solutions with a possibility of transferring the results on the multidimensional motions. As well as in the multidimensional case, the classical one-dimensional results for mono-fluids cannot be reproduced for multi-fluids automatically, in particular, due to essentially different structure of the viscous terms, namely, the presence of non-diagonal viscosity matrices; this difference in difficulty {\itshape does not depend on the dimension of the flow}. Thus, the problem under consideration is a vivid example of the role of one-dimensional theory: we temporarily put aside some difficulties (related to the dimension of the flow), and concentrate on the difficulties which do not depend on this dimension.

Well-posedness for the one-dimensional multi-fluid model with diagonal viscosity matrices was studied in \cite{mamprok.kazhpetr}, \cite{mamprok.petrov}. In the paper \cite{mamprok.nsu} we started to study the one-dimensional solvability theory for the case of the viscosity matrices of an arbitrary structure. The nearest goal is to bring the condition of the one-dimensional solvabi\-lity theory for multi-fluids to the state comparable with that for mono-fluids, i.~e. to prove the global existence of smooth solutions to the main boundary value problems.


\newpage

\begin{thebibliography}{99}
\bibitem{mamprok.feir09}  E. Feireisl, A. Novotn\'y, \emph{Singular limits in thermodynamics of viscous fluids},  Adv. Math. Fluid Mech., Birkh\"auser, Basel, 2009.

\bibitem{mamprok.kazhpetr} A.~V.~Kazhikov, A.~N.~Petrov, {\it Well-posedness of the initial-boundary value problem for a~model system of equations of a multicomponent mixture}, Din. Splosh. Sredy {\bfseries 35} (1978), 61--73. (in Russian).

\bibitem{mamprok.smj12} N.~A.~Kucher, A.~E.~Mamontov, D.~A.~Prokudin, {\it Stationary Solutions to the Equations of Dynamics of Mixtures of Heat-Conductive Compressible Viscous Fluids}, Siberian Math. J., {\bfseries 53}:6 (2012), 1075--1088.

\bibitem{mamprok.france} A.~E.~Mamontov, D.~A.~Prokudin, {\it Viscous compressible multi-fluids: modeling and multi-D existence}, Methods and Applications of Analysis, {\bfseries 20}:2 (2013), 179--195.

\bibitem{mamprok.izvran} A.~E.~Mamontov, D.~A.~Prokudin, {\it Solubility of a stationary boundary-value problem
for the equations of motion of a one-temperature mixture of viscous compressible heat--conducting fluids}, Izvestiya: Mathematics, {\bfseries 78}:3 (2014), 554--579.

\bibitem{mamprok.semi2} A.~E. Mamontov, D.~A. Prokudin, {\it Solubility of steady boundary value problem for the equations of polytropic motion of multicomponent viscous compressible fluids}, Siberian Electr. Math. Reports, {\bfseries 13} (2016), 664--693. (in Russian).

\bibitem{mamprok.semi1} A.~E. Mamontov, D.~A. Prokudin, {\it Solubility of initial boundary value problem for the equations of polytropic motion of multicomponent viscous compressible fluids}, Siberian Electr. Math. Reports, {\bfseries 13} (2016), 541--583. (in Russian).

\bibitem{mamprok.conf} A.~E. Mamontov, D.~A. Prokudin, {\it Solvability of the initial boundary value problem for the equations of viscous compressible multi-fluids},
8th International Scientific School-Conference of Young Scientists ``Theory and Numerical Methods of Solving Inverse and Ill-posed Problems'', Novosibirsk, 1-7 September, 2016.
Text of the presentation is available at
\href{http://conf.ict.nsc.ru/files/conferences/tcmiip2016/336905/abs_MamontovAE_ProkudinDA.pdf}
{{\scriptsize http://conf.ict.nsc.ru/files/conferences/tcmiip2016/336905/abs\_MamontovAE\_ProkudinDA.pdf}}

\bibitem{mamprok.smz1} A.~E. Mamontov, D.~A. Prokudin, {\it Solvability of the regularized steady problem of spatial motions of viscous compressible multi-fluids},
Siberian Math. J., {\bfseries 57}:6 (2016). To appear.

\bibitem{mamprok.smz2} A.~E. Mamontov, D.~A. Prokudin, {\it Existence of weak solutions to the problem of three-dimensional steady barotropic motions of viscous compressible multi-fluids},     Siberian Math. J., {\bfseries 58}:1 (2017).    To appear.

\bibitem{mamprok.nsu} A.~E. Mamontov, D.~A. Prokudin, {\it Local solvability of the initial boundary value problem for one-dimensional equations of polytropic flows of viscous compressible multi-fluids},    Siberian Journal of Pure and Applied Mathematics, {\bfseries 16} (2016).    To appear.

\bibitem{mamprok.izvrannew} A.~E.~Mamontov, D.~A.~Prokudin, {\it Solvability of unsteady equations of multicomponent viscous compressible fluids}, Izvestiya: Mathematics, {\bfseries 81} (2017). To appear.

\bibitem{mamprok.novstrs04} A.~Novotn\'y, I.~Stra\v skraba, {\it Introduction to the mathematical theory of compressible flow}, Oxford Lecture Series in Mathematics and Its Applications, {\bfseries 27}, Oxford University Press, Oxford, 2004.

\bibitem{mamprok.petrov} A.~N.~Petrov, {\it Well-posedness of initial-boundary value problems for one-dimensional equations of mutually penetrating flows of ideal gases}, Din. Splosh. Sredy, {\bfseries 56} (1982), 105--121. (in Russian).

\bibitem{mamprok.prokkraj} D.~A.~Prokudin, M.~V.~Krayushkina, {\it Solvability of a stationary boundary value problem for a model system of the equations of barotropic motion of a mixture of compressible viscous fluids}, Journal of Applied and Industrial Mathematics, {\bfseries 10}:3 (2016), 417--428.
\end{thebibliography}
\end{document}